\documentclass[12pt,a4paper]{amsart}
\usepackage{a4}
\usepackage[english]{babel}
\usepackage{amsmath,amssymb,amsthm}
\def\leq {\leqslant}
\def\le {\leqslant}
\def\ge {\geqslant}

\def\@bibitem[#1]#2{\item\@biblabel{#1}.\if@filesw
{\def\protect##1{\string##1\space}\immediate\write
\@auxout{\string\bibcite{#2}{#1}}}\fi\ignorespaces\@showtag{#2}}

\theoremstyle{plain}
\newtheorem{theorem}{Theorem}
\renewcommand{\theequation}%
{\arabic{section}.\arabic{equation}}
\begin{document}

\title{On $M$ -- terms approximations Besov classes in Lorentz spaces}
\author{ G. Akishev}
\address{Department of Mathematics and Information Technology, 
Buketov Karaganda State University, 
Universytetskaya~28 , 
100028, Karaganda , Republic Kazakhstan} 

\maketitle

\begin{quote}
\noindent{\bf Abstract.}
In this paper we consider Lorentz space with a mixed norm of periodic functions of many variables. We obtain the exact estimation of the best M-term approximations of Nikol'ski's, Besov's classes in the Lorentz space with the mixed norm.
\end{quote}
\vspace*{0.2 cm}

{\bf Keywords:} Lorentz space, \and Besov's class, \and approximation

 {\bf MSC:} 41A10 and 41A25

\section{Introduction}
\label{sec1}
Let  $\overline{x} =\left( x_{1} ,...,x_{m} \right) \in \Bbb{T}^{m}
=\left[ 0,2\pi \right)^{m} $  and $\theta_{j}, p_{j} \in
\left[ 1,+\infty \right)$, $j=1,...,m$. Let
$L_{\bar p, \bar{\theta}}(\Bbb{T}^{m})$ denotes the space of Lebesgue – measureable functions $ f(\bar x)$ defined on $\Bbb{R}^{m}$,
which have $ 2\pi$ -- period with respect to each variable such that
$$
\|f\|_{\overline{p},\overline{\theta}} = \|...\|f\|_{p_{1},\theta_{1}}...\|_{p_{m},\theta_{m}}
< +\infty,
$$
where
$$
\|g\|_{p,\theta} =
\left\{\int\limits_{0}^{2\pi}(g^{*}(t))^{\theta}t^{\frac{\theta}{p}-1}
dt\right\}^{\frac{1}{\theta}},
$$
where $g^{*}$ a non-increasing rearrangement of the function $|g|$ (see. [1]).

It is known that if $\theta_{j} = p_{j}, j=1,...,m$, then $L_{\bar p, \bar{\theta}}(\Bbb{T}^{m}) = L_{\bar p}(\Bbb{T}^{m})$  the Lebesgue measurable space of functions $ f(\bar x)$ defined on $\Bbb{R}^{m}$,
which have $ 2\pi$ -- period with respect to each variable with the norm
 $$
\|f\|_{\bar{p}}=\Biggl[\int_{0}^{2\pi}
\biggl[\cdots\biggl[\int_{0}^{2\pi}|f(\bar x)|
^{p_{1}}dx_{1}\biggr]^{\frac{p_{2}}{p_{1}}}
\cdots\biggr]^{\frac{p_{m}}{p_{m-1}}}dx_{m}\Biggr]^{\frac{1}{p_{m}}}<+\infty
,
 $$
where $\overline{p}=\left(p_{1},...,p_{m}\right),$  $1\leq p
_{j}<+\infty,$   $j=1,...,m$ (see [2],p. 128).

Any function $f \in L_{1}\left(\Bbb{T}^{m}\right) = L\left(\Bbb{T}^{m}\right)$ can be expanded to the Fourier series
$$
\sum\limits_{\overline{n} \in \Bbb{Z}^{m}
}a_{\overline{n}} \left( f\right) e^{i\langle\overline{n}, \overline{x}\rangle },
$$
where
$a_{\overline{n} } (f)$ Fourier coefficients of  $f\in L_{1}
\left(\Bbb{T}^{m} \right)$ with respect to multiple trigonometric system $\{e^{i\langle\overline{n}, \overline{x}\rangle}\}_{\bar n \in \Bbb{Z}^{m}},$
and $\Bbb{Z}^{m}$ is the space of points in $\Bbb{R}^{m}$ with integer coordinates.

For a function $f \in L(\Bbb{T}^{m})$ and a number $s \in \Bbb{Z}_{+} = \Bbb{N}\cup \{0\}$ let us introduce the notation
 $$
\delta_{0}(f, \bar x) = a_{0}(f), \;\;   \delta_{s} ( f,\overline{x})
=\sum\limits_{\overline{n} \in \rho (s )
}a_{\overline{n}}( f) e^{i\langle\overline{n} ,\overline{x}\rangle },
 $$
where $\langle\bar{y},\bar{x}\rangle=\sum\limits_{j=1}^{m}y_{j}
x_{j}$,
 $$
\rho (s)=\left\{ \overline{k} =( k_{1}
,...,k_{m}) \in \Bbb{Z}^{m}: \quad [2^{s -1}] \leq \max_{j=1,...,m}| k_{j}|
 <2^{s} \right\},
$$
where $[a]$ is the integer part of the number $a$.

Let us consider Nikol'skii, Besov classes( see [2], [3]). Let $1 < p_{j} < +\infty, 1 < \theta_{j} < +\infty,$ $j=1,...,m$, $1 \le \tau \le \infty$, and $r > 0$

$$
H_{\bar{p}, \bar\theta}^{r} =\left\{ f\in
L_{\bar{p}, \bar\theta}
\left(\Bbb{T}^{m} \right) : \sup\limits_{s \in \Bbb{Z}_{+}} 2^{sr} \left\|\delta_{s}(f)\right\|_{\bar{p}, \bar\theta} \leq 1\right\},
$$
 $$
B_{\bar{p}, \bar\theta, \tau}^{r} = \left\{ f\in
L_{\bar{p}, \bar\theta}
(\Bbb{T}^{m} ) :  \left(\sum\limits_{s \in \Bbb{Z}_{+}} 2^{sr\tau} \left\|\delta_{s} ( f)\right\|_{\bar{p}, \bar\theta}^{\tau} \right)^{\frac{1}{\tau}} \leq 1\right\}.
$$

It is known that for $1 \le \tau \le  \infty$ the following holds
$$
B_{\bar{p}, \bar\theta, 1}^{r} \subset B_{\bar{p}, \bar\theta, \tau}^{r}\subset  B_{\bar{p}, \bar\theta, \infty}^{r} = H_{\bar p, \bar\theta}^{r}.
$$

Let $f \in L_{\bar p, \bar\theta}(\Bbb{T}^{m})$ and $\left\{\bar{k}^{(j)} \right\}_{j=1}^{M}$ be a system of vectors $\bar{k}^{(j)} = (k_{1}^{(j)},...,k_{m}^{(j)})$ with integer coordinates. Consider the quantity
 $$
e_{M}\left(f \right)_{\bar p, \bar\theta} = \inf\limits_{\bar{k}^{(j)}, b_{j}}
\left\|f - \sum\limits_{j=1}^{M} b_{j} e^{\langle i\bar{k}^{(j)}, \bar{x}\rangle}  \right\|_{\bar p, \bar\theta},
$$
where $b_{j}$ are arbitrary numbers.
The quantity $e_{M}\left(f \right)_{\bar p, \bar\theta}$ is called the best $M$ -- term approximation of a function $f \in L_{\bar p, \bar\theta}(\Bbb{T}^{m})$.
For a given class $F \subset L_{\bar p, \bar\theta}(\Bbb{T}^{m})$ let
 $$
e_{M}\left(F \right)_{\bar p, \bar\theta} = \sup\limits_{f \in F}e_{M}\left(f \right)_{\bar p, \bar\theta}.
 $$

The best $M$ -- term  approximation was defined by S.B.Stechkin [4].
Estimations of $M$ -- term approximations of different classes were provided by R.S. Ismagilov [5], E.S. Belinsky [6], V.E. Maiorov [7], B.S. Kashin [8], R. DeVore [9], V.N. Temlyakov [10],  A.S. Romanyuk [11], Dinh Dung [12], Wang Heping and Sun Yongsheng [13], L. Q.Duan and G.S. Fang [14], W.Sickel and M. Hansen [15],  S.A. Stasyuk [16], [17]  and others (see bibliography in [18], [19], [20]).

For the case $p_{1} =... = p_{m} = p$ and $q_{1} =...=  q_{m}=\theta_{1} =...=\theta_{1} =q $ R.A. De Vore and V.N. Temlyakov [20] proved the following theorem.

{\bf Theorem A.}
\textit{Let $1 \le p, q, \tau \le \infty$ and
$r(p, q)= m\left(\frac{1}{p} - \frac{1}{q} \right)_{+}$ if $1\le p \le q\le 2$, or $1\le q \le p < \infty$ and $r(p, q) = \max \left\{\frac{m}{p}, \frac{m}{2} \right\}$ in other cases. Then for $r > r(p, q)$ the following holds
$$
e_{M}(B_{p, \tau}^{r})_{q} \asymp M^{-\frac{r}{m}+\left(\frac{1}{p} - \max \left\{\frac{1}{q}, \frac{1}{2} \right\} \right)_{+}},
$$
where $a_{+} = \max\left\{a; 0\right\}.$}

Moreover, in the case of $m(\frac{1}{p} - \frac{1}{q}) < r < \frac{m}{p}$,  $1 < p \le 2 < q < \infty$ S.A. Stasyuk [16]  proved  $e_{M}(B_{p, \tau})_{q} \asymp M^{-\frac{q}{2}(\frac{r}{m}-(\frac{1}{p} - \frac{1}{q}))}$.
In the case $r = \frac{m}{p}$  obtained $e_{M}(B_{p, \tau})_{q} \asymp M^{-\frac{1}{2}}(\log M)^{1-\frac{1}{\tau}}$ (see [17] ).

The main goal of the present paper is to find the order of the quantity
$e_{M}\left(F \right)_{\bar q, \bar\theta}$ for the class
$F = B_{\bar{p}, \bar\theta, \tau}^{r}$.

The notation $A\left( y\right) \asymp
B\left( y\right)$ means that there exist positive constants
 $C_{1},\,C_{2} $ such that  $C_{1}A\left( y\right)
\leq B\left( y\right) \leq C_{2}A\left( y\right)$. If $A \le C_{2}B$ or $A\ge C_{2}B$, then we write $A<<B$ or $A>>B$.

\section{Auxiliary results}
\label{sec2}

To  prove  the main results the following auxiliary propositions are used.

{\bf Theorem B ([21] ).} \textit{ Let $p \in (1, \infty).$ Then there exist positive numbers $C_{1}(p), C_{2}(p)$ such that for any function $f \in L_{p}(\Bbb{T}^{m})$ the following inequality holds:
$$
\|f\|_{p} << \Bigl\|\Bigl(\sum\limits_{s=0}^{\infty} |\delta_{s}(f)|^{2}\Bigr)^{\frac{1}{2}}\Bigr\|_{p} << \|f\|_{p}.
$$
}

{\bf Theorem C ([22] ).} \textit{
Let $\bar n = (n_{1},..., n_{m}), \;\;$ $n_{j} \in \Bbb{N}, j=1,...,m$ and
$$
T_{\bar n}(\bar x) = \sum\limits_{|k_{j}|\le n_{j}, \\  j=1,...,m} c_{\bar k}e^{i\langle \bar k,\bar x \rangle}.
$$
Then for $1\le p_{j} < q_{j} < \infty,$ $1 \le \theta_{j}^{(1)}, \theta_{j}^{(2)} <+\infty$, $j=1,...m$ the following inequality holds
$$
\left\|T_{\bar n} \right\|_{\bar q, \bar\theta^{(2)}}  << \prod\limits_{j=1}^{m}n_{j}^{\frac{1}{p_{j}}-\frac{1}{q_{j}}}\left\|T_{\bar n} \right\|_{\bar p, \bar\theta^{(1)}}.
$$
}

Let $\Omega_{M}$ be a set containing no more than $M$ vectors
$\bar k = (k_{1},..., k_{m})$ with integer coordinates, and
$P(\Omega_{M}, \bar x)$ be any trigonometric polynomial, which consists of harmonics with ``indices'' in $\Omega_{M}.$

{\bf Lemma 1} (see [18]). \textit{
Let $2 < q_{j}<+ \infty,$ $ j=1,...,m.$  Then for any trigonometric polynomial $P(\Omega_{N})$ and for any natural number $M < N$ there exists trigonometric polynomial $P(\Omega_{M}),$  such that the following estimation holds
$$
\|P(\Omega_{N})- P(\Omega_{M})\|_{\bar q} << (NM^{-1})^{\frac{1}{2}}\|P(\Omega_{N})\|_{2},
$$
and, moreover, $\Omega_{M}\subset\Omega_{N}$.}

\section{Main results}

Let us prove the main results.

\begin{theorem}. {\it
Let $\bar p = (p_{1},...,p_{m}),$ \quad $\bar{q} = (q_{1},...,q_{m}),$ \quad $\bar\theta^{(1)} = (\theta_{1}^{(1)},...,\theta_{m}^{(1)}),$ \quad $\bar\theta^{(2)} = (\theta_{1}^{(2)},...,\theta_{m}^{(2)})$, $1 < p_{j} \le 2 < q_{j} < \infty$, $1 < \theta_{j}^{(1)} ,\theta_{j}^{(2)} < \infty, \,\,$ $j =1,...,m, \,\,$ $1 \le \tau \le \infty$.

1. If $\sum\limits_{j=1}^{m} (\frac{1}{p_{j}}-\frac{1}{q_{j}}) <  r < \sum\limits_{j=1}^{m} \frac{1}{p_{j}}$, then
 $$
e_{M}\left(B_{\bar{p}, \bar\theta^{(1)}, \tau}^{r} \right)_{\bar q, \bar\theta^{(2)}} \asymp M^{-(2\sum\limits_{j=1}^{m} \frac{1}{q_{j}})^{-1}(r - \sum\limits_{j=1}^{m}\left(\frac{1}{p_{j}} -\frac{1}{q_{j}}\right))}. 
$$
2. If $r = \sum\limits_{j=1}^{m} \frac{1}{p_{j}}$, then
$$
e_{M}\left(B_{\bar{p},\bar\theta^{(1)} \tau}^{r} \right)_{\bar q, \bar\theta^{(2)}} \asymp M^{-\frac{1}{2}}(\log_{2} M)^{1-\frac{1}{\tau}},
$$
for $M > 1$.

3. If $r > \sum\limits_{j=1}^{m} \frac{1}{p_{j}}$, then
 $$
e_{M}\left(B_{\bar{p}, \bar\theta^{(1)}, \tau}^{r} \right)_{\bar q, \bar\theta^{(2)}} \asymp M^{-\frac{1}{m}(r + \sum\limits_{j=1}^{m}(\frac{1}{2} -\frac{1}{p_{j}}))}.
 $$
}
\end{theorem}

{\bf Proof.} Firstly, we are going to consider the upper bound in the first item. Taking into account the inclusion
$ B_{\bar{p}, \bar\theta^{(1)}, \tau}^{r} \subset H_{\bar{p}, \bar\theta^{(1)}}^{r},$ $1 \le \tau < +\infty$,
it suffices to prove it for the class $H_{\bar{p}, \bar\theta^{(1)}}^{r}.$

Let $1 < p_{j} \le 2 <  q_{j} < \infty, \;\; j = 1,...,m,$ and $\Bbb{N}$ be the set of natural numbers. For a number $M \in \Bbb{N}$ choose a natural number $n$ such that $2^{nm} < M \le 2^{(n+1)m}.$ For a function
$f \in H_{\bar{p}, \bar\theta^{(1)}}^{r}$, it is known that

$$
\|\delta_{s}(f)\|_{\bar p, \bar\theta^{(1)}} \le 2^{-sr}, \quad 1 < p_{j} < \infty, \quad j=1,...,m.
$$
We will seek an approximation polynomial $P(\Omega_{M}, \bar x)$ in the form
$$
P(\Omega_{M}, \bar x) = \sum\limits_{s=0}^{n-1}\delta_{s}(f, \bar x) +
\sum\limits_{n\le s < \alpha n}P(\Omega_{N_{s}}, \bar x),   \eqno (1)
$$
where the polynomial's $P(\Omega_{N_{s}}, \bar x)$ will be constructed for each $\delta_{s}(f, \bar x)$  in accordance with Lemma 1, and the number $\alpha > 1$ will be chosen during the construction.

Let $\sum\limits_{j=1}^{m} (\frac{1}{p_{j}}-\frac{1}{q_{j}}) <  r < \sum\limits_{j=1}^{m} \frac{1}{p_{j}}$. Suppose
$$
N_{s} = \Bigl[2^{nm}2^{s(\sum\limits_{j=1}^{m} \frac{1}{p_{j}}-r)}
2^{-n\alpha(\sum\limits_{j=1}^{m} \frac{1}{p_{j}}-r)} \Bigr] +1,
$$
where $[y]$ integer part of the number $y$.

Now we are going to show that polynomials  (1) have no more than $M$ harmonics (in terms of order).
By definition of the number $N_{s}$, we have
$$
\sum\limits_{s=0}^{n-1}\sharp \{\bar k = (k_{1},...,k_{m}) : [2^{s-1}] \le \max_{j=1,...,m}|k_{j}|
< 2^{s}\} + \sum\limits_{n \le s <\alpha n} N_{s} \le C2^{nm} +
$$
$$
+ \sum\limits_{n \le s <\alpha n}\left(2^{nm}2^{s(\sum\limits_{j=1}^{m} \frac{1}{p_{j}}-r)}
2^{-n\alpha(\sum\limits_{j=1}^{m} \frac{1}{p_{j}}-r)} +1 \right) << 2^{nm} +
(\alpha  - 1)n << 2^{nm} \asymp M,
$$
where $\sharp A$ denotes the number of elements of the set $A$.

Next, by the property of the norm
we have
$$
\|f - P(\Omega_{M})\|_{\bar q, \bar\theta^{(2)}} \le  \left\|\sum\limits_{n \le s <\alpha n} (\delta_{s}(f) - P(\Omega_{N_{s}})) \right\|_{\bar q, \bar\theta^{(2)}}
+
$$
$$
+ \left\|\sum\limits_{\alpha n \le s < +\infty} \delta_{s}(f) \right\|_{\bar q, \bar\theta^{(2)}}=J_{1}(n) + J_{2}(n). \eqno (2)
$$

Let us estimate $J_{2}(n).$ Applying the inequality of different metrics for trigonometric polynomials (see Theorem C), we can  obtained
$$
J_{2}(n) \le \sum\limits_{\alpha n \le s < +\infty} \|\delta_{s}(f)\|_{\bar q, \bar\theta^{(2)}}
<< \sum\limits_{\alpha n \le s < +\infty}2^{s\sum\limits_{j=1}^{m}(\frac{1}{p_{j}} -\frac{1}{q_{j}})} \|\delta_{s}(f)\|_{\bar p, \bar\theta^{(1)}}.
$$
Therefore, taking into account $f \in H_{\bar p, \bar\theta^{(1)}}^{r}$ and   $\sum\limits_{j=1}^{m}(\frac{1}{p_{j}} -\frac{1}{q_{j}}) < r$, we get
$$
J_{2}(n) << \sum\limits_{\alpha n \le s < +\infty}2^{-s(r - \sum\limits_{j=1}^{m}(\frac{1}{p_{j}} -\frac{1}{q_{j}}))} << 2^{-n\alpha (r - \sum\limits_{j=1}^{m}(\frac{1}{p_{j}} -\frac{1}{q_{j}}))}.  \eqno (3)
$$
Let us estimate $J_{1}(n)$.
Using the property of the quasi-norm, Lemma 1 and the inequality of different metrics (see Theorem C), we get
 $$
J_{1}(n) = \left\|\sum\limits_{n \le s <\alpha n} (\delta_{s}(f) - P(\Omega_{N_{s}})) \right\|_{\bar q, \bar\theta^{(2)}} << \sum\limits_{n \le s <\alpha n} \left\|\delta_{s}(f) - P(\Omega_{N_{s}})\right\|_{\bar q, \bar\theta^{(2)}} <<
$$
$$
<< \sum\limits_{n \le s <\alpha n} (N_{s}^{-1}2^{sm})^{\frac{1}{2}}\left\|\delta_{s}(f)\right\|_{2} <<
 $$
 $$
<< \sum\limits_{n \le s <\alpha n} (N_{s}^{-1}2^{sm})^{\frac{1}{2}}2^{s\sum\limits_{j=1}^{m}(\frac{1}{p_{j}} -\frac{1}{2})}\left\|\delta_{s}(f)\right\|_{\bar p, \bar\theta^{(1)}} <<
$$
$$
<<  \sum\limits_{n \le s <\alpha n} N_{s}^{-\frac{1}{2}}2^{s\sum\limits_{j=1}^{m}\frac{1}{p_{j}}}2^{-sr} <<
$$
$$
<< 2^{-\frac{nm}{2}}2^{\frac{n\alpha}{2}(\sum\limits_{j=1}^{m}\frac{1}{p_{j}}-r )}
\sum\limits_{n \le s <\alpha n}2^{s(\sum\limits_{j=1}^{m}\frac{1}{p_{j}}-r)\frac{1}{2}} <<
2^{-\frac{nm}{2}}2^{n\alpha(\sum\limits_{j=1}^{m}\frac{1}{p_{j}}-r )}.
\eqno (4)
$$
Suppose $\alpha = m(2\sum\limits_{j=1}^{m}\frac{1}{q_{j}})^{-1}$. Then from the inequality (4), we get
$$
J_{1}(n) \le C 2^{-nm(2\sum\limits_{j=1}^{m}\frac{1}{q_{j}})^{-1}(r - \sum\limits_{j=1}^{m}(\frac{1}{p_{j}} -\frac{1}{q_{j}}))} \asymp M^{-(2\sum\limits_{j=1}^{m}\frac{1}{q_{j}})^{-1}(r - \sum\limits_{j=1}^{m}(\frac{1}{p_{j}} -\frac{1}{q_{j}}))}. \eqno (5)
$$
For $\alpha = m(2\sum\limits_{j=1}^{m}\frac{1}{q_{j}})^{-1}$, using the inequality (3) and taking into account $2^{nm} \asymp M$ we obtain
$$
J_{2}(n) << M^{-(2\sum\limits_{j=1}^{m}\frac{1}{q_{j}})^{-1}(r - \sum\limits_{j=1}^{m}(\frac{1}{p_{j}} -\frac{1}{q_{j}}))}. \eqno (6)
$$
By (5) and (6), we get from the inequality (2) the following
$$
\|f - P(\Omega_{M})\|_{\bar q, \bar\theta^{(2)}} << M^{-(2\sum\limits_{j=1}^{m}\frac{1}{q_{j}})^{-1}(r - \sum\limits_{j=1}^{m}(\frac{1}{p_{j}} -\frac{1}{q_{j}}))}
$$
for any function $f \in H_{\bar p, \bar\theta^{(1)}}^{r}$ in the case of $\sum\limits_{j=1}^{m} (\frac{1}{p_{j}}-\frac{1}{q_{j}}) <  r < \sum\limits_{j=1}^{m} \frac{1}{p_{j}}$.

From the inclusion $B_{\bar p, \bar\theta^{(1)}, \tau}^{r} \subset H_{\bar p, \bar\theta^{(1)}}^{r}$ and the definition of $M$ -- term approximation, it follows that
$$
e_{M}\left(B_{\bar{p}, \bar\theta^{(1)}, \tau}^{r} \right)_{\bar q, \bar\theta^{(2)}} <<
M^{-(2\sum\limits_{j=1}^{m}\frac{1}{q_{j}})^{-1}(r - \sum\limits_{j=1}^{m}(\frac{1}{p_{j}} -\frac{1}{q_{j}}))}
$$
in the case of $\sum\limits_{j=1}^{m} (\frac{1}{p_{j}}-\frac{1}{q_{j}}) <  r < \sum\limits_{j=1}^{m} \frac{1}{p_{j}}$.

Let us consider the lower bound. We will use the well-known formula ([23], p.25)
$$
e_{M}(f)_{\bar q, \bar\theta^{(2)}} = \inf\limits_{\Omega_{M}} \sup\limits_{P\in L^{\perp}, \\ \|P\|_{\bar{q}', \bar{\theta^{(2)}}'} \le 1}
\left|\int_{\Bbb{T}^{m}} f(\bar x)P(\bar x)d\bar{x}\right|, \eqno (7)
$$
where $\bar{q}^{'} = (q_{1}^{'},...,q_{m}^{'}),\;\;$ $\bar{\theta}^{(2)^{'}} = (\theta_{1}^{(2)^{'}},...,\theta_{m}^{(2)^{'}}),\;\;$ $\frac{1}{q_{j}} + \frac{1}{q_{j}^{'}} = 1, \;\;$ $\frac{1}{\theta_{j}^{(2)}} + \frac{1}{\theta_{j^{(2)^{'}}}} = 1$, $j=1,...,m,$ and $L_{M}^{\perp}$ is the set of functions that are orthogonal to the subspace of  trigonometric polynomials with harmonics in the set $\Omega_{M}$.

Consider the function
$$
F_{\bar q, n}(\bar x) = \sum\limits_{\max\limits_{j=1,...,m} |k_{j}|\le 2^{[nm(2\sum\limits_{j=1}^{m}\frac{1}{q_{j}})^{-1}]}}
e^{i\langle \bar{k}, \bar{x} \rangle}.
$$
Let $\Omega_{M}$ be a set of $M$ vectors with integer coordinates. Suppose
$$
g(\bar x) = F_{\bar q, n}(\bar x) -  \sum\limits_{\bar k \in \Omega_{M}}^{*}
e^{i\langle \bar{k}, \bar{x} \rangle},
$$
where the sum $\sum\limits_{\bar k \in \Omega_{M}}^{*}$ contains those terms in the function $F_{\bar q, n}(\bar x)$ with indices only in $\Omega_{M}$. By the inequality
$$
\|\sum\limits_{\max\limits_{j=1,...,m} |k_{j}|\le 2^{l}} e^{i\langle \bar{k}, \bar{x} \rangle}\|_{\bar p, \bar\theta^{(1)}} << 2^{l\sum\limits_{j=1}^{m}(1-\frac{1}{p_{j}})} \eqno (8)
$$
and the Perseval's equality for $1 < q_{j}^{'} < 2, j=1,...,m$, we obtain
$$
\|g\|_{\bar{q}^{'}, \bar{\theta}^{(2)^{'}}} \le \|F_{\bar q, n}\|_{\bar{q}^{'}, \bar{\theta}^{(2)^{'}}} + (2\pi)^{\sum\limits_{j=1}^{m}(\frac{1}{q_{j}} - \frac{1}{2})}\|\sum\limits_{\bar k \in \Omega_{M}}^{*}e^{i\langle \bar{k}, \bar{x} \rangle} \|_{2} <<
$$
$$
<< (2^{\frac{nm}{2}} + M^{\frac{1}{2}}) << 2^{\frac{nm}{2}}.  \eqno (9)
$$
Now we consider the function
$$
P_{1}(\bar x) = C_{2}2^{-\frac{nm}{2}}g(\bar x).    \eqno (10)
$$
Then  (9) implies  follows that the function $P_{1}$ satisfies the assumptions of the formula (7) for some constant $C_{2} > 0$.

Consider the function
$$
f_{n}(\bar x) = C_{3}2^{-nm(2\sum\limits_{j=1}^{m}\frac{1}{q_{j}})^{-1}(r - \sum\limits_{j=1}^{m}(\frac{1}{p_{j}} -1))}F_{\bar q, n}(\bar x).
$$
By the inequality (8), we get
 $$
\sum\limits_{s=0}^{\infty} 2^{sr}\left\|\delta_{s}(f_{n})\right\|_{\bar{p}, \bar{\theta}^{(1)}} <<
$$
$$
 << 2^{-nm(2\sum\limits_{j=1}^{m}\frac{1}{q_{j}})^{-1}(r - \sum\limits_{j=1}^{m}(\frac{1}{p_{j}} -1))}\sum\limits_{s=0}^{[m(2\sum\limits_{j=1}^{m}\frac{1}{q_{j}})^{-1}]} 2^{sr}\left\|\delta_{s}(F_{\bar q, n})\right\|_{\bar{p}, \bar{\theta}^{(1)}} <<
$$
$$
 << 2^{-nm(2\sum\limits_{j=1}^{m}\frac{1}{q_{j}})^{-1}(r - \sum\limits_{j=1}^{m}(\frac{1}{p_{j}} -1))}\sum\limits_{s=0}^{[m(2\sum\limits_{j=1}^{m}\frac{1}{q_{j}})^{-1}]} 2^{sr}
2^{\sum\limits_{j=1}^{m}(1 - \frac{1}{p_{j}})} \le C_{3}.
$$
Hence, the function $C_{3}^{-1}f_{n} \in B_{\bar p, \bar{\theta}^{(1)}, 1}^{r}$.

For the functions (10) and (11), we have by the formula (7), the following
$$
e_{M}(f_{n})_{\bar q, \bar{\theta}^{(2)}} >> \inf\limits_{\Omega_{M}}
\left|\int_{\Bbb{T}^{m}} f_{n}(\bar x)P_{1}(\bar x)d\bar{x}\right| >>
$$
$$
>> 2^{-nm(2\sum\limits_{j=1}^{m}\frac{1}{q_{j}})^{-1}(r - \sum\limits_{j=1}^{m}(\frac{1}{p_{j}} -1))} 2^{-\frac{nm}{2}}(\|F_{\bar q, n}\|_{2}^{2} - M) >>
$$
$$
>> 2^{-nm(2\sum\limits_{j=1}^{m}\frac{1}{q_{j}})^{-1}(r - \sum\limits_{j=1}^{m}(\frac{1}{p_{j}} -1))} 2^{-\frac{nm}{2}}2^{nm(2\sum\limits_{j=1}^{m}\frac{1}{q_{j}})^{-1}} =
$$
$$
= C 2^{-nm(2\sum\limits_{j=1}^{m}\frac{1}{q_{j}})^{-1}(r - \sum\limits_{j=1}^{m}(\frac{1}{p_{j}} -\frac{1}{q_{j}}))}. \eqno (12)
$$
Hence,  it follows from (12) by the inclusion $B_{\bar p, \bar{\theta}^{(1)}, 1}^{r} \subset B_{\bar p, \bar{\theta}^{(1)}, \tau}^{r}$ that
$$
e_{M}(f_{n})_{\bar q, \bar{\theta}^{(2)}} >>
2^{-nm(2\sum\limits_{j=1}^{m}\frac{1}{q_{j}})^{-1}(r - \sum\limits_{j=1}^{m}(\frac{1}{p_{j}} -\frac{1}{q_{j}}))}
$$
in the case of $\sum\limits_{j=1}^{m} (\frac{1}{p_{j}}-\frac{1}{q_{j}}) <  r < \sum\limits_{j=1}^{m} \frac{1}{p_{j}}$.

So we have proved the first item.

Now we consider the case $r = \sum\limits_{j=1}^{m} \frac{1}{p_{j}}$.
Let $f \in B_{\bar p, \bar{\theta}^{(1)}, \tau}^{r}$. Suppose \\$\alpha = m(2\sum\limits_{j=1}^{m}\frac{1}{q_{j}})^{-1}$ and
$$
N_{s} = \Bigl[2^{nm}n^{\frac{1}{\tau} - 1}\left\|\delta_{s}(f_{n})\right\|_{\bar{p}, \bar{\theta}^{(1)}} 2^{sr}\Bigr] +1.
$$
Then, by definition of the numbers $N_{s}$ and by Holder's inequality, we obtain
$$
\sum\limits_{s=0}^{n-1}\sharp \rho (s) + \sum\limits_{n \le s <\alpha n} N_{s} <<
$$
$$
<< 2^{nm} + (\alpha -1)n +  2^{nm}n^{\frac{1}{\tau} - 1}\sum\limits_{n \le s <\alpha n}\left\|\delta_{s}(f_{n})\right\|_{\bar{p}, \bar{\theta}^{(1)}} 2^{sr} <<
$$
$$
<< 2^{nm} + (\alpha -1)n +  2^{nm}n^{\frac{1}{\tau} - 1}((\alpha -1)n)^{\frac{1}{\tau'}}\left(\sum\limits_{s =0}^{\infty}\left\|\delta_{s}(f_{n})\right\|_{\bar{p}, \bar{\theta}^{(1)}}^{\tau} 2^{sr\tau}  \right)^{\frac{1}{\tau}} <<
$$
$$
<< 2^{nm} \asymp M.
$$

To estimate $J_{1}(n)$ let  $\beta = \max\{q_{1},...,q_{m}\}.$ Then $\beta >2$ and $L_{\beta}(\Bbb{T}^{m}) \subset L_{\bar{q}, \bar\theta^{(2)}}(\Bbb{T}^{m}).$ Therefore by applying Theorem B and by the norm property we obtain
$$
J_{1}(n) = \left\|\sum\limits_{n \le s <\alpha n} (\delta_{s}(f) - P(\Omega_{N_{s}})) \right\|_{\bar{q}, \bar\theta^{(2)}} \le C\left\|\sum\limits_{n \le s <\alpha n} (\delta_{s}(f) - P(\Omega_{N_{s}}))\right\|_{\beta} <<
$$
$$
<< \left\|\left(\sum\limits_{n \le s <\alpha n} |\delta_{s}(f) - P(\Omega_{N_{s}})|^{2} \right)^{\frac{1}{2}}\right\|_{\beta} <<
$$
$$
<< \left(\sum\limits_{n \le s <\alpha n} \left\|\delta_{s}(f) - P(\Omega_{N_{s}})\right\|_{\beta}^{2} \right)^{\frac{1}{2}}.
$$
It implies by Lemma 1 and by the inequality of different metrics (see Theorem C) that
$$
J_{1}(n) << \left(\sum\limits_{n \le s <\alpha n} N_{s}^{-1}2^{sm}\left\|\delta_{s}(f)\right\|_{2}^{2} \right)^{\frac{1}{2}} <<
$$
$$
<< \left(\sum\limits_{n \le s <\alpha n} N_{s}^{-1}2^{sm}2^{2s\sum\limits_{j=1}^{m}(\frac{1}{p_{j}} -\frac{1}{2})}\|\delta_{s}(f)\|_{\bar{p}, \bar\theta^{(1)}}^{2} \right)^{\frac{1}{2}}.
$$

Next, since $r = \sum\limits_{j=1}^{m} \frac{1}{p_{j}}$, we have  by definition of the numbers $N_{s}$ and using Holder's inequality, the following
$$
J_{1}(n) << \left(\sum\limits_{n \le s <\alpha n} N_{s}^{-1}2^{sm}2^{2s\sum\limits_{j=1}^{m}(\frac{1}{p_{j}} -\frac{1}{2})}\left\|\delta_{s}(f)\right\|_{\bar p, \bar{\theta}^{(1)}}^{2} \right)^{\frac{1}{2}} <<
$$
$$
<< (2^{-nm}n^{1 - \frac{1}{\theta}})^{\frac{1}{2}}
\left(\sum\limits_{n \le s <\alpha n}2^{sr}\left\|\delta_{s}(f)\right\|_{\bar p, \bar{\theta}^{(1)}}
\right)^{\frac{1}{2}} <<
$$
$$
<< (2^{-nm}n^{1 - \frac{1}{\theta}})^{\frac{1}{2}}
\left(\sum\limits_{n \le s <\alpha n}2^{sr\theta}\left\|\delta_{s}(f)\right\|_{\bar p, \bar{\theta}^{(1)}}^{\tau}
\right)^{\frac{1}{2\tau}} (\alpha - 1)^{\frac{1}{2\tau'}} =
$$
$$
= C 2^{-\frac{nm}{2}}n^{1 - \frac{1}{\tau}} \asymp M^{-\frac{1}{2}}(\log M)^{1 - \frac{1}{\tau}}.
$$
Thus
$$
J_{1}(n) << M^{-\frac{1}{2}}(\log M)^{1 - \frac{1}{\tau}}  \eqno (13)
$$
in the case of $r = \sum\limits_{j=1}^{m} \frac{1}{p_{j}}$.

For the estimation of $J_{2}(n)$ we apply Holder's inequality and taking into account that $r = \sum\limits_{j=1}^{m} \frac{1}{p_{j}}$ and $\alpha = m(2\sum\limits_{j=1}^{m}\frac{1}{q_{j}})^{-1}$,
we obtain
$$
J_{2}(n) << \sum\limits_{n\alpha \le s < +\infty}2^{s\sum\limits_{j=1}^{m} (\frac{1}{p_{j}} - \frac{1}{q_{j}})}\left\|\delta_{s}(f)\right\|_{\bar p, \bar{\theta}^{(1)}} <<
$$
$$
<< \left(\sum\limits_{s = 0}^{\infty}2^{sr\tau}\left\|\delta_{s}(f)\right\|_{\bar p, \bar{\theta}^{(1)}}^{\tau}
\right)^{\frac{1}{\theta}} \left(\sum\limits_{n\alpha \le s < +\infty} 2^{-s\tau'\sum\limits_{j=1}^{m} \frac{1}{q_{j}}}
\right)^{\frac{1}{\tau'}} <<
$$
$$
<< 2^{-n\alpha\sum\limits_{j=1}^{m} \frac{1}{q_{j}}} = C2^{-\frac{nm}{2}} \asymp M^{-\frac{1}{2}}, \eqno (14)
$$
where $\tau' = \frac{\tau}{\tau - 1}$.

By $(13)$ and $(14)$ the inequality (2) implies that
$$
\|f - P(\Omega_{M})\|_{\bar q, \bar{\theta}^{(2)}} << M^{-\frac{1}{2}}(\log M)^{1 - \frac{1}{\tau}}
$$
in the case $r = \sum\limits_{j=1}^{m} \frac{1}{p_{j}}$. It proves the upper bound estimation in the second item.

Let  $r > \sum\limits_{j=1}^{m} \frac{1}{p_{j}}$. Suppose
$$
N_{s} = \Bigl[2^{n(r - \sum\limits_{j=1}^{m} (\frac{1}{p_{j}} -1))}  2^{-s(r - \sum\limits_{j=1}^{m} \frac{1}{p_{j}})}\Bigr] +1.
$$
Then
$$
\sum\limits_{s=0}^{n-1}\sharp \rho(s) + \sum\limits_{n \le s <\alpha n} N_{s} <<
$$
$$
<< 2^{nm} + (\alpha -1)n + 2^{n(r - \sum\limits_{j=1}^{m} (\frac{1}{p_{j}} -1))}\sum\limits_{n \le s <\alpha n} 2^{-s(r - \sum\limits_{j=1}^{m} \frac{1}{p_{j}})} <<
$$
$$
<< 2^{nm} + (\alpha -1)n << 2^{nm} << M.
$$
If $f \in H_{\bar p, \bar{\theta}^{(1)}}^{r}$, then, by definition of the numbers $N_{s}$ and $r > \sum\limits_{j=1}^{m} \frac{1}{p_{j}}$
we obtain
$$
J_{1}(n) \le \left(\sum\limits_{n \le s <\alpha n} N_{s}^{-1}2^{sm}2^{2s\sum\limits_{j=1}^{m}(\frac{1}{p_{j}} -\frac{1}{2})}\left\|\delta_{s}(f)\right\|_{\bar p, \bar{\theta}^{(1)}}^{2} \right)^{\frac{1}{2}} <<
$$
$$
<< 2^{-\frac{n}{2}(r - \sum\limits_{j=1}^{m} (\frac{1}{p_{j}} -1))}
\left(\sum\limits_{n \le s <\alpha n} 2^{s(r+\sum\limits_{j=1}^{m}\frac{1}{p_{j}})}\left\|\delta_{s}(f)\right\|_{\bar p, \bar{\theta}^{(1)}}^{2} \right)^{\frac{1}{2}} <<
$$
$$
<< 2^{-\frac{n}{2}(r - \sum\limits_{j=1}^{m} (\frac{1}{p_{j}} -1))}
\left(\sum\limits_{n \le s <\alpha n} 2^{-s(r-\sum\limits_{j=1}^{m}\frac{1}{p_{j}})} \right)^{\frac{1}{2}} <<
2^{-n(r+\sum\limits_{j=1}^{m}(\frac{1}{2} - \frac{1}{p_{j}}))}.
$$
Thus,
$$
J_{1}(n) << M^{-\frac{1}{m}(r+\sum\limits_{j=1}^{m}(\frac{1}{2} - \frac{1}{p_{j}}))} \eqno(15)
$$
in the case of $r > \sum\limits_{j=1}^{m} \frac{1}{p_{j}}$.

To estimate $J_{2}(n)$, we suppose $\alpha = (r+\sum\limits_{j=1}^{m}(\frac{1}{2} - \frac{1}{p_{j}}))(r+\sum\limits_{j=1}^{m}(\frac{1}{q_{j}} - \frac{1}{p_{j}}))^{-1} $ and get
$$
J_{2}(n) << \sum\limits_{n\alpha \le s < \infty} 2^{s\sum\limits_{j=1}^{m}(\frac{1}{p_{j}} - \frac{1}{q_{j}})}\left\|\delta_{s}(f)\right\|_{\bar p, \bar{\theta}^{(1)}} \le \sum\limits_{n\alpha \le s < \infty} 2^{-s(r+ \sum\limits_{j=1}^{m}(\frac{1}{q_{j}} - \frac{1}{p_{j}}))} <<
$$
$$
<< 2^{-n\alpha(r+ \sum\limits_{j=1}^{m}(\frac{1}{q_{j}} - \frac{1}{p_{j}}))}
<< 2^{-n(r+ \sum\limits_{j=1}^{m}(\frac{1}{2} - \frac{1}{p_{j}}))} << M^{-\frac{1}{m}(r+ \sum\limits_{j=1}^{m}(\frac{1}{2} - \frac{1}{p_{j}}))}
\eqno (16)
$$
for a function $f \in H_{\bar p, \bar{\theta}^{(1)}}^{r}$. By (15) and (14), it follows from (2) that
$$
\|f - P(\Omega_{M})\|_{\bar q, \bar{\theta}^{(2)}} << M^{-\frac{1}{m}(r+ \sum\limits_{j=1}^{m}(\frac{1}{2} - \frac{1}{p_{j}}))}
$$
for any function $f \in H_{\bar p, \bar{\theta}^{(1)}}^{r}$ in the case of
$r > \sum\limits_{j=1}^{m} \frac{1}{p_{j}}$.

From $B_{\bar p, \bar{\theta}^{(1)}, \tau}^{r} \subset H_{\bar p, \bar{\theta}^{(1)}}^{r}$ it follows that
$$
e_{M}\left(B_{\bar{p}, \bar{\theta}^{(1)}, \tau}^{r} \right)_{\bar q, \bar{\theta}^{(2)}} << M^{-\frac{1}{m}(r+ \sum\limits_{j=1}^{m}(\frac{1}{2} - \frac{1}{p_{j}}))}
 $$
in the case of
$r > \sum\limits_{j=1}^{m} \frac{1}{p_{j}}$. It proves the upper bound estimation in the item 3.

Let us consider the lower bound estimation in the case $r = \sum\limits_{j=1}^{m} \frac{1}{p_{j}}$.
Consider the function
$$
g_{1}(\bar x) = \sum\limits_{s=1}^{n}\sum\limits_{\bar k \in \rho(s)} \prod\limits_{j=1}^{m} k_{j}^{-1}\cos k_{j}x_{j}. \eqno (17)
$$
Then
$$
\delta_{s}(g_{1}, \bar x) = \sum\limits_{\bar k \in \rho(s)} \prod\limits_{j=1}^{m} k_{j}^{-1}\cos k_{j}x_{j}.
$$
It is known that for a function $d_{s}(\bar x) = \sum\limits_{\bar k \in \rho(s)} \prod\limits_{j=1}^{m}\cos k_{j}x_{j}$ the following relation holds
$$
\|d_{s}\|_{\bar p, \bar{\theta}^{(1)}} \asymp 2^{s\sum\limits_{j=1}^{m}(1 - \frac{1}{p_{j}})},
\;\; 1 < p_{j}, \theta_{j}^{(1)} < +\infty , \;\; j=1,...,m.
$$
Therefore by the inequality of distinct metrics (see Theorem C) and by the Marcinkiewicz's theorem on multipliers, we have
$$
\|\delta_{s}(g_{1})\|_{\bar p, \bar{\theta}^{(1)}} << 2^{-sm}\|d_{s}\|_{\bar p, \bar{\theta}^{(1)}} \le C
2^{-s\sum\limits_{j=1}^{m}\frac{1}{p_{j}}}.
$$
Hence, since $r = \sum\limits_{j=1}^{m} \frac{1}{p_{j}}$
we obtain
$$
\left(\sum\limits_{s=0}^{\infty} 2^{sr\tau}\left\|\delta_{s}(g_{1})\right\|_{\bar p, \bar{\theta}^{(1)}}^{\tau} \right)^{\frac{1}{\tau}} \le C_{1}n^{\frac{1}{\tau}}.
$$
Therefore the function $f_{1}(\bar x) = C_{1}^{-1}n^{-\frac{1}{\tau}}g_{1}(\bar x)$
belongs to the class $B_{\bar{p}, \bar{\theta}^{(1)}, \tau}^{r}, \;\; 1< p_{j} < +\infty, j=1,...,m.$

Now, we are going to construct a function $P_{1}$, which satisfies the conditions of the formula (7).
Let
$$
v_{1}(\bar x) = \sum\limits_{s=1}^{n}\sum\limits_{\bar k \in \rho(s)} \prod\limits_{j=1}^{m}\cos k_{j}x_{j}
$$
and $\Omega_{M}$ be an arbitrary set of vectors $\bar k = (k_{1},...,k_{m})$ in $M$ with integer coordinates. Consider the function
 $$
u_{1}(\bar x) = \sum\limits_{s=1}^{n}\sum\limits_{\bar k \in \rho(s)\cap\Omega_{M}} \prod\limits_{j=1}^{m}\cos k_{j}x_{j} .
 $$
Suppose $w_{1}(\bar x) = v_{1}(\bar x) - u_{1}(\bar x)$. Then, since $1 < q_{j}^{'} = \frac{q_{j}}{q_{j} - 1} < 2,\;\; j=1,...,m$, we obtain,  by the Perseval`s equality, the following
$$
\|w_{1}\|_{\bar{q}^{'}, \bar{\theta}^{(2)^{'}}} \le \|v_{1}\|_{\bar{q}^{'}, \bar{\theta}^{(2)^{'}}} + \|u_{1}\|_{2} \le \|v_{1}\|_{\bar{q}^{'}, \bar{\theta}^{(2)^{'}}} + CM^{\frac{1}{2}}.
$$
By the property of quasi-norm and the estimation of the norm of the Dirichlet kernel in the Lorentz space, we have
$$
\|v_{1}\|_{\bar{q}^{'}, \bar{\theta}^{(2)^{'}}} << \sum\limits_{s=1}^{n} \left\|\delta_{s}(g_{1})\right\|_{\bar{q}^{'}, \bar{\theta}^{(2)^{'}}} <<
$$
$$
<< \sum\limits_{s=1}^{n} 2^{s\sum\limits_{j=1}^{m}(1 - \frac{1}{q_{j}^{'}})}
<< 2^{n\sum\limits_{j=1}^{m}(1 - \frac{1}{q_{j}^{'}})} = C2^{n\sum\limits_{j=1}^{m}\frac{1}{q_{j}}}.
$$
Therefore, taking into account $\frac{1}{q_{j}} < \frac{1}{2}, \;\; j =1,...,m$, we get
$$
\|w_{1}\|_{\bar{q}^{'}, \bar{\theta}^{(2)^{'}}} << (2^{\frac{nm}{2}} + M^{\frac{1}{2}}) \le C_{2}2^{\frac{nm}{2}}.
$$
Hence the function
$$
P_{1}(\bar x) = C_{2}^{-1}2^{-\frac{nm}{2}}w_{1}(\bar x)
$$
satisfies the conditions of the formula (7). Then, by substituting the functions $f_{1} and  P_{1}$ into
(7) and by orthogonally of a trigonometric system, we obtain
$$
e_{M}\left(f_{1}\right)_{\bar q, \bar{\theta}^{(2)}} >> \sum\limits_{n_{1} \le s < n}
\sum\limits_{\bar k \in \rho(s)} \prod\limits_{j=1}^{m} k_{j}^{-1}
2^{-\frac{nm}{2}}n^{-\frac{1}{\tau}} >>
$$
$$
\ge C (\ln 2)^{m}\sum\limits_{n_{1} \le s < n} 2^{-\frac{nm}{2}}n^{-\frac{1}{\tau}} = C (\ln 2)^{m}2^{-\frac{nm}{2}}n^{-\frac{1}{\tau}}(n - n_{1}) \ge
$$
$$
>> (\ln 2)^{m}2^{-\frac{nm}{2}}n^{1 - \frac{1}{\tau}} \asymp M^{-\frac{1}{2}}(\log_{2} M)^{1 - \frac{1}{\tau}},
$$
where  $n_{1}$ is a natural number such that $n_{1} < \frac{n}{2}$.

So, for the function $f_{1} \in B_{\bar{p}, \bar{\theta}^{(1)}, \tau}^{r}$ it has been proved that
$$
e_{M}\left(f_{1}\right)_{\bar q, \bar{\theta}^{(2)}} >> M^{-\frac{1}{2}}(\log_{2} M)^{1 - \frac{1}{\tau}}
$$
in the case of $r = \sum\limits_{j=1}^{m} \frac{1}{p_{j}}$.
Hence
$$
e_{M}\left(B_{\bar{p}, \bar{\theta}^{(1)}, \tau}^{r}\right)_{\bar q, \bar{\theta^{(2)}}} >> M^{-\frac{1}{2}}(\log_{2} M)^{1 - \frac{1}{\tau}}
$$
in the case of $r = \sum\limits_{j=1}^{m} \frac{1}{p_{j}}$.
It proves the lower bound estimation in the second item.

Let us prove the lower bound estimation for the case $r > \sum\limits_{j=1}^{m} \frac{1}{p_{j}}$.
Since in this case an upper bound estimation of the quantity $e_{M}\left(B_{\bar{p}, \bar{\theta}^{(1)}, \tau}^{r}\right)_{\bar q, \bar{\theta}^{(2)}}$ does not depend on $\tau$ and $B_{\bar{p}, \bar{\theta}^{(1)}, 1}^{r} \subset B_{\bar{p}, \bar{\theta}^{(1)},\tau}^{r}$, $1 < \tau < + \infty$, it suffices to prove the lower bound estimation for $B_{\bar{p}, \bar{\theta}^{(1)}, 1}^{r}$.

For a number $M \in \Bbb{N}$, we choose a natural number $n$ such that $2^{nm} < M \le 2^{(n+1)m}$ and $2M \le \sharp\rho(n)$, where $\sharp\rho(n)$ denotes the number of elements in the set $\rho(n)$.

Consider the following function
$$
f_{3}(\bar x) =n^{-1} \sum\limits_{s=1}^{n} 2^{-s\sum\limits_{j=1}^{m}(1 - \frac{1}{p_{j}})} \sum\limits_{\bar k \in \rho(s)}  \prod\limits_{j=1}^{m} k_{j}^{-\frac{r}{m}}\cos k_{j}x_{j}.
$$
Then
$$
\|\delta_{s}(f_{3})\|_{\bar p, \bar{\theta}^{(1)}} << 2^{-sr}n^{-1}.
$$
Hence
$$
\sum\limits_{s=0}^{\infty} 2^{sr}\left\|\delta_{s}(f_{3})\right\|_{\bar p, \bar{\theta}^{(1)}}
\le C_{3}
$$
i.e. the function $C_{3}^{-1}f_{3} \in B_{\bar{p}, \bar{\theta}^{(1)}, 1}^{r}$.

Next, consider the functions
$$
v_{3}(\bar x) = \sum\limits_{s=1}^{n}\sum\limits_{\bar k \in \rho(s)} \prod\limits_{j=1}^{m}\cos k_{j}x_{j},
$$
$$
u_{3}(\bar x) = \sum\limits_{s=1}^{n}\sum\limits_{\bar k \in \rho(s)\cap\Omega_{M}} \prod\limits_{j=1}^{m}\cos k_{j}x_{j}.
$$
Suppose $w_{3}(\bar x) = v_{3}(\bar x) - u_{3}(\bar x).$ By the Perseval`s equality,
$$
\|u_{3}\|_{2} \le M^{\frac{1}{2}},
$$
$$
\|v_{3}\|_{2} = 2^{\frac{(n-1)m}{2}}.
$$
From these relations, we obtain, by the properties of the norm, the following
$$
\|w_{3}\|_{2} \le \|v_{3}\|_{2} + \|u_{3}\|_{2} \le C_{4}2^{\frac{nm}{2}}.
$$
Therefore the function $P_{3}(\bar x) = C_{4}^{-1}2^{-\frac{nm}{2}}w_{3}(\bar x)$ satisfies the conditions of the formula (7). Since $2 < q_{j} \;\; j=1,...,m$, we have
$e_{M}\left(f_{3}\right)_{2} \le Ce_{M}\left(f_{3}\right)_{\bar q, \bar{\theta}^{(2)}}$.
Now, by the formula (7), we get
$$
e_{M}\left(f_{3}\right)_{\bar q, \bar{\theta}^{(2)}} >> e_{M}\left(f_{3}\right)_{2} >> n^{-1}
2^{-\frac{nm}{2}}\sum\limits_{s=1}^{n} 2^{-s\sum\limits_{j=1}^{m}(1 - \frac{1}{p_{j}})} \sum\limits_{\bar k \in \rho(s)}  \prod\limits_{j=1}^{m} k_{j}^{-\frac{r}{m}} >>
$$
$$
>> n^{-1}2^{-\frac{nm}{2}}\sum\limits_{s=1}^{n} 2^{-s\sum\limits_{j=1}^{m} (1 - \frac{1}{p_{j}})}2^{s(m-r)} =
$$
$$
= C n^{-1}
2^{-\frac{nm}{2}}\sum\limits_{s=1}^{n} 2^{-s(r -\sum\limits_{j=1}^{m} \frac{1}{p_{j}})} >> 2^{-n(r + \sum\limits_{j=1}^{m}(\frac{1}{2} - \frac{1}{p_{j}}))}.
$$
It follows from the relation $2^{nm} \asymp M$ that
$$
e_{M}\left(f_{3}\right)_{\bar q, \bar{\theta}^{(2)}} >> M^{-\frac{1}{m}(r + \sum\limits_{j=1}^{m}(\frac{1}{2} - \frac{1}{p_{j}})) }
$$
in the case $r > \sum\limits_{j=1}^{m} \frac{1}{p_{j}}$ for the function
$C_{3}^{-1}f_{3} \in B_{\bar{p}, \bar{\theta^{(1)}}, 1}^{r}$.
Hence
$$
e_{M}\left(B_{\bar{p}, \bar{\theta^{(1)}}, 1}^{r}\right)_{\bar q, \bar{\theta}^{(2)}} >> M^{-\frac{1}{m}(r + \sum\limits_{j=1}^{m}(\frac{1}{2} - \frac{1}{p_{j}})) }.
$$
Therefore
$$
e_{M}\left(B_{\bar{p}, \bar{\theta}^{(1)}, \tau}^{r}\right)_{\bar q} >> M^{-\frac{1}{m}(r + \sum\limits_{j=1}^{m}(\frac{1}{2} - \frac{1}{p_{j}})) }
$$
in the case $r > \sum\limits_{j=1}^{m} \frac{1}{p_{j}}$.
So Theorem 1 has been proved.

\begin{theorem}.
{\it Let $\bar p = (p_{1},...,p_{m})$, $\bar q = (q_{1},...,q_{m}),$ $1 < p_{j} < q_{j} \le 2$, $1 < \theta_{j}^{(1)} ,\theta_{j}^{(2)} < \infty, \,\,$ $j =1,...,m, \,\,$ $1\le \tau \le +\infty$.

If $r > \sum\limits_{j=1}^{m}(\frac{1}{p_{j}} - \frac{1}{q_{j}}),$ then
$$
e_{M}\left(B_{\bar{p}, \bar\theta^{(1)}, \tau}^{r}\right)_{\bar q, \theta^{(2)}} \asymp
M^{-\frac{1}{m}(r - \sum\limits_{j=1}^{m}(\frac{1}{p_{j}} - \frac{1}{q_{j}}))}.
$$}
\end{theorem}

{\bf Proof.} For a number $M \in \Bbb{N}$ choose a natural number $n$ such that $M \asymp 2^{nm}$. By the inequality of distinct metrics and by Holder`s inequality, we have
$$
\|f - \sum\limits_{s=0}^{n} \delta_{s}(f)\|_{\bar q, \bar\theta^{(2)}} \le
\sum\limits_{s=n}^{\infty} \|\delta_{s}(f)\|_{\bar q, \bar\theta^{(2)}} \le
$$
$$
\le \Bigl[\sum\limits_{s=0}^{\infty} 2^{s\tau r}\|\delta_{s}(f)\|_{\bar q, \bar\theta^{(2)}}^{\tau} \Bigr]^{\frac{1}{\tau}} \le \Bigl[   \sum\limits_{s=n}^{\infty} 2^{s\tau'(r - \sum\limits_{j=1}^{m}(\frac{1}{p_{j}} - \frac{1}{q_{j}}))}\Bigr]^{\frac{1}{\tau'}} <<
$$
$$
<< 2^{n(r - \sum\limits_{j=1}^{m}(\frac{1}{p_{j}} - \frac{1}{q_{j}}))} \le C M ^{-\frac{1}{m}(r - \sum\limits_{j=1}^{m}(\frac{1}{p_{j}} - \frac{1}{q_{j}}))}
$$
for $f \in B_{\bar{p}, \bar\theta^{(1)}, \tau}^{r}, \;\; \frac{1}{\tau} + \frac{1}{\tau^{'}} = 1.$
Therefore
$$
e_{M}(f)_{\bar q, \bar\theta^{(2)}} \le
\|f - \sum\limits_{s=0}^{n} \delta_{s}(f)\|_{\bar q, \bar\theta^{(2)}} <<
M ^{-\frac{1}{m}(r - \sum\limits_{j=1}^{m}(\frac{1}{p_{j}} - \frac{1}{q_{j}}))}.
$$
Hence
$$
e_{M}\left(B_{\bar{p}, \bar\theta^{(1)}, \tau}^{r}\right)_{\bar q, \bar\theta^{(2)}} <<
M ^{-\frac{1}{m}(r - \sum\limits_{j=1}^{m}(\frac{1}{p_{j}} - \frac{1}{q_{j}}))}.
$$
It proves the upper bound estimation.

For the lower bound estimation, let us consider the function
$$
f_{0}(\bar x) = n^{- r + \sum\limits_{j=1}^{m}(\frac{1}{p_{j}} - 1)}V_{n}(\bar x),
$$
where $V_{n}(\bar x)$ is a Valle-Poisson sum with multiplicity.

Next,  following the proof in [19] (pp. 46-47) and applying Theorem B, we obtain the lower bound estimation of the quantity $e_{M}\left(B_{\bar{p}, \bar\theta^{(1)}, \tau}^{r}\right)_{\bar q, \bar{\theta}^{(2)}}$.

\begin{theorem}.   
{\it Let $\bar p = (p_{1},...,p_{m})$, $\bar q = (q_{1},...,q_{m}),$ $2  \le p_{j} < q_{j} < \infty $,  $1 < \theta_{j}^{(1)} ,\theta_{j}^{(2)} < \infty, \,\,$ $j =1,...,m, \,\,$ $1\le \tau \le +\infty$.
If $r >\frac{m}{2},$ then
$$
e_{M}\left(B_{\bar{p}, \bar\theta^{(1)}, \tau}^{r}\right)_{\bar q, \bar\theta^{(2)}} \asymp
M^{-\frac{r}{m}}.
$$}
\end{theorem}
{\bf Proof.} By the inclusion
$B_{\bar{p}, \bar\theta^{(1)}, \tau}^{r} \subset B_{\bar{2}, \bar\theta^{(1)}, \tau}^{r} \subset  H_{2, \bar\theta^{(1)}}^{r}$,
we have
$$
e_{M}\left(B_{\bar{p}, \bar\theta^{(1)}, \tau}^{r}\right)_{\bar q, \bar\theta^{(2)}} \le e_{M}\left(B_{2, \bar\theta^{(1)}, \tau}^{r}\right)_{\bar q, \bar\theta^{(2)}} \le e_{M}\left(H_{2, \bar\theta^{(1)}}^{r}\right)_{\bar q, \bar\theta^{(2)}}.
$$
By Theorem 1,
$$
e_{M}\left(H_{2, \bar\theta^{(1)}}^{r}\right)_{\bar q, \bar\theta^{(2)}} << M^{-\frac{r}{m}}.
$$
for $p_{j} = 2, j=1,...,m.$ Hence
$$
e_{M}\left(B_{\bar{p}, \bar\theta^{(1)}, \tau}^{r}\right)_{\bar q, \bar\theta^{(2)}} << M^{-\frac{r}{m}}.
$$
it proves the upper bound estimation.

Let us consider the lower bound estimation. Consider Rudin-Shapiro's polynomial (see [24], p.155) of the type
$$
R_{s}(x) = \sum\limits_{s=2^{s-1}}^{2^{s}} \varepsilon_{k}e^{ikx}, \;\; x \in [0, 2\pi], \;\; \varepsilon_{k} = \pm 1.
$$
it is known that $\|R_{s}\|_{\infty} = \max\limits_{x\in [0, 2\pi]}|R_{s}(x)| << 2^{\frac{s}{2}}$ (see [24], p. 155). For a given number $M$ choose a number $n$ such that $M \asymp 2^{nm}.$ Now consider the function
$$
f_{0}(\bar x) = 2^{-n(\frac{m}{2} + r)}\sum\limits_{s=1}^{n}\prod_{j=1}^{m} R_{s}(x_{j})
$$
Then,by the continuity, $f_{0} \in L_{\bar p, \bar\theta^{(1)}}(\Bbb{T}^{m})$ and
$$
\sum\limits_{s=0}^{\infty} 2^{s\tau r}\|\delta_{s}(f_{0})\|_{\bar p, \bar\theta^{(1)}}^{\tau}
= 2^{-n(\frac{m}{2} + r)}\sum\limits_{s=1}^{n}2^{s\tau r}\|\prod\limits_{j=1}^{m} R_{s}(x_{j})\|_{\bar p, \bar\theta^{(1)}}^{\tau} \le
$$
$$
\le 2^{-n(\frac{m}{2} + r)}\sum\limits_{s=1}^{n}2^{s(\frac{m}{2} + r)\tau} \le C_{0}.
$$
Hence, the function $C_{0}^{-1}f_{0} \in B_{\bar{p}, \bar\theta^{(1)}, \tau}^{r}$. Now construct a function $P(\bar x)$, which would satisfy the conditions in the formula (7). Suppose
$$
v_{0}(\bar x) = \sum\limits_{s=1}^{n}\prod\limits_{j=1}^{m} R_{s}(x_{j}), \quad
u_{0}(\bar x) = \sum\limits_{s}^{*}\prod\limits_{j=1}^{m} R_{s}(x_{j}),
 $$
where the sign $*$ means that the polynomial $u_{0}(\bar x)$ contains only those harmonics of $v_{0}$ which have indices in $\Omega_{M}$. Suppose $w_{0}(\bar x) = v_{0}(\bar x) - u_{0}(\bar x).$ Then, since $1 < q_{j}^{'} =\frac{q_{j}}{q_{j}-1} < 2, \;\; j=1,...,m$, and by the Perceval's equality, we have
$$
\|w_{0}\|_{\bar{q}^{'}, \bar{\theta}^{(2)^{'}}} \le \|w_{0}\|_{2} \le C_{1}2^{\frac{nm}{2}}.
$$
Therefore, for the function $P_{0}(\bar x) = C_{1}^{-1}2^{-\frac{nm}{2}}w_{0}(\bar x)$, the inequality holds $\|P_{0}\|_{\bar{q}^{'}, \bar{\theta^{(2)^{'}}}} \le 1.$ Now using the formula (7), we obtain
$$
e_{M}\left(B_{\bar{p}, \bar\theta^{(1)} \tau}^{r}\right)_{\bar q, \bar\theta^{(2)}} >> e_{M}(f_{0})_{\bar q, \bar\theta^{(2)}}
>> 2^{-n(\frac{m}{2} + r)}2^{-\frac{nm}{2}}(2^{nm} - M) >>
$$
$$
>> 2^{-n(m + r)}2^{nm} >> M^{-\frac{r}{m}}.
$$
So
$$
e_{M}\left(B_{\bar{p},\bar\theta^{(1)}, \tau}^{r}\right)_{\bar q, \bar\theta^{(2)}} >> M^{-\frac{r}{m}}.
$$
It proves Theorem 3.

{\bf Corollary.}
Let $1 < p \le 2 < q < \infty$, $1\le \tau \le \infty$ and $r = \frac{m}{p}$. Then
$$
e_{M}\left(B_{p, \tau}^{r} \right)_{q} \asymp M^{-\frac{1}{2}}(\log M)^{1-\frac{1}{\theta}}.
$$

The proof follows from the second item of Theorem 2.1 if $p_{j} = \theta_{j}^{(1)} = p,$ $q_{j} = \theta_{j}^{(2)} = q, j = 1,...,m.$

{\bf Remark.} In the case $p_{j} = \theta_{j}^{(1)} = p,$ $q_{j} = \theta_{j}^{(2)} = q, j = 1,...,m$ and $r > m(\frac{1}{p} - \frac{1}{q}),$, the results of R.A. DeVore and V.N. Temlyakov [20] follow from Theorem 1 - 3. If $1 < p \le 2 < q < \infty$ and $ m(\frac{1}{p} - \frac{1}{q}) < r \le \frac{m}{p},$ the results of S.A. Stasyuk [16], [17] follow from the first and second items of Theorem 1.

The cases $p_{j} = \theta_{j}^{(1)},$ $q_{j} = \theta_{j}^{(2)}, j = 1,...,m.$ of Theorem 1 - 3 were announced in [25] and in of Theorem 1 the first item proved [26].

\end{document}